\documentclass[a4paper,11pt]{amsart}
\usepackage{float}
\usepackage{amssymb}
\usepackage{mathrsfs}
\usepackage{amsmath,amssymb,amsthm,latexsym,amscd,mathrsfs}
\usepackage{indentfirst}
\usepackage{stmaryrd}
\usepackage{graphicx}
\usepackage{extarrows}
\usepackage{multirow}
\usepackage{latexsym}
\usepackage{amsfonts}
\usepackage{color}
\usepackage{pictexwd,dcpic}
\usepackage{graphicx}
\usepackage{psfrag}
\usepackage{hyperref}
\usepackage{comment}
\usepackage{enumerate}

\usepackage{amsmath}
\usepackage{amsfonts}
\usepackage{amssymb}
\usepackage{amsthm}
\usepackage{mathrsfs}
\usepackage{indentfirst}
\usepackage{bbm}
\usepackage{bm}

\numberwithin{equation}{section} \theoremstyle{plain}
\newtheorem{thm}{Theorem}[section]

\newtheorem*{thm EdMc}{Edholm-McNeal Theorem}

\newtheorem{lem}[thm]{Lemma}
\newtheorem{cor}[thm]{Corollary}
\newtheorem{rem}[thm]{Remark}

\makeatletter

\def\<{\langle}
\def\>{\rangle}
\def\({\left(}
\def\){\right)}
\def\[{\left[}
\def\]{\right]}

\makeatother

	\title[$L^p$ boundedness of the Bergman projection]{$L^p$ boundedness of the Bergman projection on generalized Hartogs triangle in $\mathbb{C}^{n+1}$}

\author[Q. Fu]{Qian Fu$^{*}$}
\address{$^{*}$Teaching Experiment Platform, Beijing Normal University, Zhuhai 519087, P. R. China.}
\email{qianfu@bnu.edu.cn}

\author[G. T. Deng]{Guantie Deng${ }^{\dagger}$}
\address{${ }^{\dagger}$School of Mathematical Sciences, Beijing Normal University, Beijing 100875, P.R. China.}\email{96022@bnu.edu.cn}

\thanks{${ }^{\dagger}$Corresponding author.
}

\allowdisplaybreaks
\begin{document}
	
\maketitle

\begin{abstract}
In this paper, we study generalized Hartogs triangle of exponent $\gamma>0$, $\Omega^{n+1}_\gamma =\{(z,w)\in \mathbb{C}^n\times \mathbb{C}: |z|^\gamma < |w| < 1\}$, and 
obtain a sharp range of $p$ for the boundedness of the Bergman projection on the domain considered here. It generalizes the results of Edholm and McNeal [J. Geom. Anal. \textbf{27}, 2658-2683 (2017)] for $n = 1$ to any dimension $n$.\\

\noindent{\bf Key words}:\ Hartogs triangle, $L^p$ regularity, Bergman kernel, Bergman projection.\\
{\bf MSC 2020}:\ 32A36;\ 32A25.
\end{abstract}

\section{Introduction}
 \setlength{\parindent}{2em} Let $\Omega$ be a domain in $\mathbb{C}^{n}$ and $\mathcal{O}(\Omega)$ be the space of holomorphic functions on $\Omega$.  For $p>0$, denote
$
L^{p}(\Omega)=\left\{f:\left(\int_{\Omega}|f|^{p} d V\right)^{\frac{1}{p}}:=\|f\|_{p}<\infty\right\},
$
where $dV(z)$ is the ordinary Lebesgue volume measure on $\Omega$. 
For $p=2$, $L^{2}(\Omega)$ is a Hilbert space with the inner product:
\begin{equation}\label{innerproduct}
\langle f,g\rangle=\int_{\Omega}f(z)\overline{g(z)}dV(z).
\end{equation}
Let $A^{p}(\Omega)=\mathcal{O}(\Omega) \cap L^{p}(\Omega)$. It follows from Bergman inequality that $A^{2}(\Omega)$ is a closed subspace of $L^{2}(\Omega).$ The Bergman projection associated to $\Omega$ is the orthogonal projection
 $\mathbf{P}_{\Omega}: L^{2}(\Omega) \longrightarrow A^{2}(\Omega)$, which has an integral representation
$$
\mathbf{P}_{\Omega} f(z)=\int_{\Omega} B_{\Omega}(z, w) f(w) d V(w), 
$$
for all $f \in L^{2}(\Omega) $ and $z \in \Omega$.
Here the function $B_{\Omega}(z, w)$ defined on $\Omega\times\Omega$ is the Bergman kernel.
For more discussion of the Bergman kernel, please see \cite{CD,DG2,KZ}.
The Bergman projection is a linear operator by definition, and
 it is self-adjoint with respect to the inner product (\ref{innerproduct}). See also \cite{SGK2} for further definitions and basic properties.



Different types of regularity of the Bergman projection are of particular interest. In general, the regularity of $\mathbf{P}_{\Omega}$ depends closely on the geometry of $\Omega$. 
For various geometric conditions on $\Omega$, understanding the range of $p$ for which $\mathbf{P}_{\Omega}$ is $L^p$ bounded is an active area of research. 
Smoothly bounded domains with various convexity conditions on the boundary were considered, see for example \cite{MCN,DHP}.
In \cite{YE}, the author constructed pseudoconvex domains in $\mathbb{C}^2$, where $\mathbf{P}_{\Omega}$ is bounded if and only if $p = 2$. Beberok \cite{TB} also considered the $L^p$ boundedness of the Bergman projection on the following generalization of the Hartogs triangle:
$\mathcal{H}^{n+1}_k := \{(z, w) \in \mathbb{C}^n\times \mathbb{C}: ||z|| < |w|^k < 1\},$
where $k \in \mathbb{Z}^+$ and $||\cdot||$ is the Euclidean norm in $\mathbb{C}^n$. On some other domains, the projection has only a finite range of mapping regularity (cf. \cite{CL1,LD2}, etc).

The Hartogs triangle $\left\{\left(z_{1}, z_{2}\right) \in \mathbb{C}^{2}:\left|z_{1}\right|<\left|z_{2}\right|<1\right\}$ is a pseudoconvex domain that is the source of many counterexamples in several complex variables; see \cite{SMC}.
The generalized Hartogs triangles recently studied by Edholm and McNeal \cite{LD3,LD1} are a class of pseudoconvex domains in $\mathbb{C}^{2}$ defined for $\gamma>0$ by
$
\mathbb{H}_{\gamma}=\left\{\left(z_{1}, z_{2}\right) \in \mathbb{C}^{2}:\left|z_{1}\right|^{\gamma}<\left|z_{2}\right|<1\right\}$.
$\mathbb{H}_{\gamma}$ exhibits the same pathological behavior as the classical Hartogs triangle due to the singularity of the boundary, which is non-Lipschitz at the origin, with the additional surprising dependence on the rationality or irrationality of the power $\gamma$. More specifically,
\begin{thm EdMc}
	\label{thm Edholm and McNeal Bergman projection}
	The Bergman projection $\mathbf{P}_{\gamma}\ (\gamma\in\mathbb{R^+})$ is a bounded operator from $L^{p}(\mathbb{H}_{\gamma})$ to $A^{p}(\mathbb{H}_{\gamma})$ if and only if $p\in\mathscr{A}$.
	\begin{enumerate}[(1)]
		\item If $\gamma=m / l$ is rational, then $\mathscr{A}=\left(\frac{2 m+2 l}{m+l+1}, \frac{2 m+2 l}{m+l-1}\right)$
		$\left( m, l \in \mathbb{Z}^{+},\ \operatorname{gcd} (m, l)=1\right)$;
		\item If $\gamma$ is irrational, then
		$\mathscr{A}=\{2\}$.
	\end{enumerate}
\end{thm EdMc}
Recently, Chen et al. \cite{CL3}  proved that if $\Omega$ can be covered by the polydisc  through a surjective proper rational holomorphic map, then $P_\Omega$ is $L^p$-bounded for $p\in (r,r')$, where $r<2$ and $r'>2$ are conjugate exponents. 
This result can be applied to 
certain generalized Hartogs triangles $\mathbb{H}_{\gamma} \,(\gamma \in\mathbb{Q}^+)$ to conclude that the Bergman projection of these domains is $L^p$-regular for $p$ in a finite interval.


In this article, we mainly study the following bounded regions. For $\gamma >0$, we define the domain $\Omega^{n+1}_\gamma\subseteq \mathbb{C}^n\times \mathbb{C}$ by
$$\Omega^{n+1}_\gamma =\{(z,w)\in \mathbb{C}^n\times \mathbb{C}: |z|^\gamma < |w| < 1\}$$
and call $\Omega^{n+1}_\gamma$ the generalized Hartogs triangle of exponent $\gamma,n$. On the generalized Hartogs triangle, we denote the Bergman projection by $\mathbf{P}_{\gamma}$ and the Bergman kernel as $B_{\gamma,n}((z,w),(s,t))$. As usual, the operator 
$\mathbf{P}_{\gamma}$ is extended to supersets of  $L^2(\Omega^{n+1}_\gamma)$ by setting 
$$
\mathbf{P}_{\gamma} (f)(z,w)=\int_{\Omega^{n+1}_\gamma} B_{\gamma,n}((z,w),(s,t)) f(s,t) d V(s,t), 
$$
whenever the integral is defined. The primary purpose of this paper is to show that the Bergman projection of $\Omega^{n+1}_\gamma$, $\mathbf{P}_{\gamma}$, is $L^p$ bounded for only a restricted range of $p \in (1, \infty).$ The precise statement of our main result that extends \cite{LD1} for $n = 1$ to any dimension $n$ is as follows.
\begin{thm}[\bf Main Theorem]
	\label{thm Bergman projection}
	The Bergman projection $\mathbf{P}_{\gamma}\ (\gamma\in\mathbb{R^+})$ is a bounded operator from $L^{p}(\Omega^{n+1}_\gamma)$ to $A^{p}(\Omega^{n+1}_\gamma)$ if and only if $p\in\mathscr{A}$.
	\begin{enumerate}[(1)]
		\item If $\gamma=m / l$ is rational, then $\mathscr{A}=\left(\frac{2 m+2 nl}{m+nl+1}, \frac{2 m+2 nl}{m+nl-1}\right)$
		$\left( m, l \in \mathbb{Z}^{+},\ \operatorname{gcd} (m, l)=1\right)$;
		\item If $\gamma$ is irrational, then
		$\mathscr{A}=\{2\}$.
	\end{enumerate}
\end{thm}



\subsection*{Outline of the proof}
We will prove the main theorem above in three steps, corresponding to Sections \ref{Qboundedness}, \ref{Qnonboundedness}, \ref{Irrational}, respectively:

\begin{enumerate}[]
\item \textbf{Step 1.} We divide $\gamma$ into rational and irrational numbers.
 In the former case, we decompose the Bergman space $A^{2}\left(\Omega^{n+1}_\gamma \right)$ to obtain an estimate of the Bergman kernel (see Theorem \ref{thm3}).  We then use Schur's Lemma in Section \ref{Qboundedness} to discuss the boundedness of the Bergman projection operator,  thus obtaining the range of $p$ in Theorem \ref{thm Bergman projection} (1) (see Theorem \ref{necessity}).
 \item \textbf{Step 2.}
 For the proof of unboundedness, in Section \ref{Qnonboundedness} we exhibit a single function $f \in L^{\infty}(\Omega^{n+1}_\gamma)$ such that $\mathbf{P}_\gamma f \notin L^p(\Omega^{n+1}_\gamma)$.
 The range of $L^p$ boundedness in Theorem \ref{thm Bergman projection} (1) is shown to be sharp (see Theorem \ref{Sufficiency}). Hence, this completes the proof of Theorem   \ref{thm Bergman projection}  (1).  
 \item \textbf{Step 3.}
 When $\gamma$ is irrational, the Bergman kernel of $B_{\gamma,n}((z,w),(s,t))$ is not a rational function, so in Section \ref{Irrational} we prove  Theorem \ref{thm Bergman projection} (2) using Dirichlet's theorem on rational approximation of $\gamma \notin \mathbb{Q}$. 
\end{enumerate}




\section{The Rational Case: $L^{p}$ Boundedness}\label{Qboundedness}
\subsection{Bergman Space $A^{2}\left(\Omega^{n+1}_\gamma \right)$ Decomposition}
The following lemmas are needed before decomposing the Bergman Space.
\begin{lem}{\rm\cite{HB}}\label{lemmaHB}
	\label{formula}
	For any $v_1,\cdots, v_n \geq 0$,
	$$
	\int_{S^{2n-1}}|\zeta_1|^{2v_1}\cdots|\zeta_n|^{2v_n}d\sigma(\zeta)=\frac{2v!\pi^n}{\Gamma(n+|v|)},
	$$
	where $|v| = v_1 +\cdots+ v_n$, $v! = \Gamma(v_1 + 1)\cdots \Gamma(v_n + 1),$
	and $S^{2n-1}$ is  the unit sphere in $\mathbb{C}^n$ with respect to the surface measure $d\sigma$.
\end{lem}

\begin{lem}\label{A2}
	If the monomial $z^\alpha w^\beta \in A^2(\Omega^{n+1}_\gamma),$ 
	then multi-indices $(\alpha,\beta)\in \mathbb{N}^n\times \mathbb{Z}$ should meet $|\alpha| + \gamma (\beta + 1) > -n$, and 
	$$
	\|z^\alpha w^\beta\|^2_{L^2(\Omega^{n+1}_\gamma)} =\frac{\pi^{n+1}\alpha!}{ (\beta+\frac{|\alpha|+n}{\gamma}+1)\Gamma(|\alpha|+n+1)}.
	$$
\begin{proof}
	Using polar coordinates and Lemma \ref{lemmaHB}, we obtain
    \begin{align}\label{L2}
	\|z^\alpha w^\beta\|^2_{L^2(\Omega^{n+1}_\gamma)}
	&=\int_{\Omega^{n+1}_\gamma}|z^\alpha|^{2}|w|^{2\beta}dV(z,w)\notag\\
	&= \int_{0<|w|<1}|w|^{2\beta}\int_{0}^{|w|^{\frac{1}{\gamma}}}\int_{|\zeta|=1}|\zeta_1|^{2\alpha_1}\cdots|\zeta_n|^{2\alpha_n}d\sigma(\zeta)r^{2|\alpha|+2n-1}drdV(w)\notag\\
	&= \frac{\pi^n \alpha!}{(|\alpha|+n)\Gamma(|\alpha|+n)}\int_{0<|w|<1}|w|^{2\beta+\frac{2|\alpha|+2n}{\gamma}}dV(w)\notag\\
	&= \frac{2\pi^{n+1}\alpha!}{ \Gamma(|\alpha|+n+1)}\int_0^1 r^{2\beta+\frac{2|\alpha|+2n}{\gamma}+1}dr.
    \end{align}	
This integral converges if and only if $2\beta+\frac{2|\alpha|+2n}{\gamma}+1>-1$, $i.e.$, $ |\alpha| + \gamma (\beta + 1) > -n$. Furthermore, when the integral (\ref{L2}) converges, it equals
$$
\frac{2\pi^{n+1}\alpha!}{ \Gamma(|\alpha|+n+1)} \cdot \frac{1}{2\beta+\frac{2|\alpha|+2n}{\gamma}+2}
=\frac{\pi^{n+1}\alpha!}{ (\beta+\frac{|\alpha|+n}{\gamma}+1)\Gamma(|\alpha|+n+1)}.
$$
\end{proof}	
\end{lem}
Now, we consider the Bergman kernel $B_{\gamma,n}((z,w),(s,t))$, $z,s \in\mathbb{C}^n, w,t \in\mathbb{C}$. We will use the biholomorphic transformation law of the Bergman kernel (see, e.g., \cite[Proposition 1.4.12]{SGK2}). Let $\Omega_1, \Omega_2$ be domains in $\mathbb{C}^n$. Let $f: \Omega_1 \rightarrow \Omega_2$ be biholomorphic. Then
$$
\operatorname{det} J_{\mathbb{C}} f(z) K_{\Omega_2}(f(z), f(\zeta)) \operatorname{det} \overline{J_{\mathbb{C}} f(\zeta)}=K_{\Omega_1}(z, \zeta),
$$ 
where $K_{\Omega}$ denotes the Bergman kernel of the domain $\Omega$, $J_{\mathbb{C}} f(z)$ denotes the complex Jacobian matrix of the biholomorphic mapping $f$. 
In particular, unitary transformations preserve the Bergman kernel. Choose one unitary matrix $U$  such that $z=|z|\mathbbm{1}U^{-1},$ where $\mathbbm{1}=(1,0,\cdots,0)$.  Then, from Lemma \ref{A2}, we obtain
\begin{align}\label{kernel}
	B_{\gamma,n}((z,w),(s,t))
	&=\sum\limits_{|\alpha| + \gamma (\beta + 1) > -n}\frac{z^{\alpha}w^{\beta}(\overline{s^\alpha t^\beta})}{\|z^\alpha w^\beta\|^2_{L^2(\Omega^{n+1}_\gamma)}}
	=\sum\limits_{|\alpha| + \gamma (\beta + 1) > -n}\frac{(|z|\mathbbm{1})^{\alpha}w^{\beta}(\overline{(sU)^\alpha t^\beta})}{\|z^\alpha w^\beta\|^2_{L^2(\Omega^{n+1}_\gamma)}}\notag\\
	&=\sum\limits_{(\alpha_1,\beta)\in \Lambda_{\gamma,n}}\frac{(|z|\mathbbm{1}\overline{U}^{T}\overline{s}^{T})^{\alpha_1}(w\cdot \overline{t})^{\beta}}{N_{\gamma,n}(\alpha_1,\beta)}
	=\sum\limits_{(\alpha_1,\beta)\in \Lambda_{\gamma,n}}\frac{(z\cdot\overline{s})^{\alpha_1}(w\cdot\overline{t})^{\beta}}{N_{\gamma,n}(\alpha_1,\beta)},
\end{align}
where
\begin{align}\label{Lambda1}
	\Lambda_{\gamma,n}=\{(\alpha_1,\beta)\in \mathbb{N}\times \mathbb{Z} :\alpha_1+\gamma(\beta+1)>-n\},
\end{align}  
\begin{align}\label{norm}
N_{\gamma,n}(\alpha_1,\beta)=\int_{\Omega^{n+1}_\gamma}|z_1|^{2\alpha_1}|w|^{2\beta}dV(z,w)=\frac{\pi^{n+1}\Gamma(\alpha_1+1)}{ (\beta+\frac{\alpha_1+n}{\gamma}+1)\Gamma(\alpha_1+n+1)}.
\end{align} 
Therefore, below, we only need to consider a two-dimensional array that satisfies the inequality (\ref{Lambda1}).

When $\gamma = \frac{m}{l}\in \mathbb{Q}^+,$ $gcd(m,l)=1,$ the strict inequality defining $(\alpha_1,\beta) \in \Lambda_{\gamma,n}$ can be re-expressed as a non-strict inequality:
\begin{align}\label{Lambda2}
	\Lambda_{\frac{m}{l},n}
	&=\{(\alpha_1,\beta)\in \mathbb{N}\times \mathbb{Z} :\alpha_1+\frac{m}{l}(\beta+1)>-n\}\notag\\
	&=\{(\alpha_1,\beta)\in \mathbb{N}\times \mathbb{Z} :l\alpha_1+m\beta\geq-m-nl+1\}.
\end{align} 
We split the Bergman space into $m$ orthogonal subspaces
\begin{align}\label{subspaces}
A^{2}\left(\Omega^{n+1}_{m/l} \right)=\mathcal{S}_{0} \oplus \mathcal{S}_{1} \oplus \cdots \oplus \mathcal{S}_{m-1},
\end{align}
where $\mathcal{S}_{j}$ is the subspace spanned by monomials of the form $z^{\alpha}w^{\beta}\in A^{2}\left(\Omega^{n+1}_{m/l} \right)$, where $\alpha_{1} = j \bmod m$. Let
\begin{align}\label{Gj}
\mathcal{G}_{j}=\left\{(\alpha_1,\beta) \in \Lambda_{\frac{m}{l},n}: \alpha_{1} = j \bmod m\right\},
\end{align}
 $\mathcal{G}_{j} \cap \mathcal{G}_{k}=\varnothing$ if $j \neq k$. Each $\mathcal{S}_{j}$ is a closed subspace of $A^{2}\left(\Omega^{n+1}_{m/l} \right)$,  and thus a Hilbert space. Therefore the orthogonal projection $\mathcal{K}_{j},$ $L^{2}\left(\Omega^{n+1}_{m/l}\right) \longrightarrow \mathcal{S}_{j}$, is well defined and represented by integration against a kernel, $K_{j}$. It follows that
\begin{align}\label{decompose}
B_{\frac{m}{l},n}((z,w),(s,t))
=\sum\limits_{j=0}^{m-1}K_{j}((z,w),(s,t)).
\end{align}
Call each $K_{j}$ a sub-Bergman kernel. 
In the next subsection, we shall focus on the subspaces $\mathcal{S}_{j}$ and estimate each $K_{j}$. 
We use the following notation to simplify writing various inequalities. If $A$ and
$B$ are functions depending on several variables, we write $A\lesssim B$ to signify that there exists a
constant $K > 0$, independent of relevant variables, such that $A \leqslant K B$. The independence
of which variables will be clear in context. We also write $A\thickapprox B$ to mean that $A\lesssim B \lesssim A$. If $x \in \mathbb{R}$, $\lfloor x \rfloor$ will denote the greatest integer not exceeding $x$.

\subsection{Estimation of the Bergman Kernel} 
Let $\gamma=\frac{m}{l} \in \mathbb{Q}^{+},$
$ \operatorname{gcd}(m, l)=1$. For each $j=0, \ldots, m-1$, let $K_{j}$ be the sub-Bergman kernel of $B_{\frac{m}{l},n}$. By definition,
$\{z^{\alpha}w^{\beta}: (\alpha_1,\beta) \in \mathcal{G}_{j}\}$ is an orthonormal basis for $\mathcal{S}_{j}$, where $\mathcal{G}_{j}$ is given by (\ref{Gj}).
Then, (\ref{kernel}) implies that $K_{j}$ can be written as the following sum:
\begin{align}\label{Kj}
	K_{j}((z,w),(s,t))
	=\sum\limits_{(\alpha_1,\beta)\in \mathcal{G}_{j}}\frac{(z\cdot\overline{s})^{\alpha_1}(w\cdot\overline{t})^{\beta}}{N_{\frac{m}{l},n}(\alpha_1,\beta)},
\end{align}	
where  $N_{\frac{m}{l},n}(\alpha_1,\beta)$ by (\ref{norm}).
\begin{thm}
	\label{thm3}
Let $m, l \in \mathbb{Z}^{+}$ be relatively prime. The sub-Bergman kernel $K_{j}$ of domain $\Omega^{n+1}_{m/l}$ satisfies the estimate
$$
|K_{j}((z,w),(s,t))|\lesssim \frac{|b|^{\frac{lj}{m}+(n+1)l-E_j-1}}{|1-b|^{2}|b^l-a^{m}|^{n+1}},
$$
where 
$a=z\cdot\overline{s},$ 
$ b=w\cdot\overline{t},$ 
$ E_{j}=\left\lfloor\frac{(j+n)l-1}{m}\right\rfloor$.
\end{thm}

\begin{proof} First we find $K_{j}((z,w),(z,w))$ and then use polarization to move off the diagonal. Let $a=|z|^{2}, b=|w|^{2}$. 
Starting from (\ref{Kj}) and using (\ref{Lambda2}),
\begin{align}\label{Kj1}
	&K_{j}((z,w),(z,w))\notag\\
    =&\sum\limits_{(\alpha_1,\beta)\in \mathcal{G}_{j}}\frac{a^{\alpha_1}b^{\beta}}{N_{\frac{m}{l},n}(\alpha_1,\beta)}\notag \\
	=&\frac{1}{\pi^{n+1}} \sum_{\alpha_{1} \in \mathcal{R}_{j}} \sum_{\beta \in \mathcal{N}_{j}}\left[\frac{\Gamma(\alpha_1+n+1)(\alpha_1+n)l}{\Gamma(\alpha_1+1)m}+\frac{\Gamma(\alpha_1+n+1)(\beta+1)}{\Gamma(\alpha_1+1)}\right] a^{\alpha_{1}} b^{\beta},
\end{align}
where $\mathcal{R}_{j}:=\left\{\alpha_{1} \geq 0: \alpha_{1}=j \bmod m\right\}$ and 
$\mathcal{N}_{j}:=\left\{\beta \in \mathbb{Z}:\beta \geq-\frac{l \alpha_{1}}{m}+\frac{1-nl-m}{m}\right\}$. Let $\ell(j)$ represent the smallest integer not less than $-\frac{l \alpha_{1}}{m}+\frac{1-nl-m}{m}$.
Notice that
$$
-\frac{l \alpha_{1}}{m}+\frac{1-nl-m}{m}=-1-\frac{l\left(\alpha_{1}-j\right)}{m}-\frac{(j+n)l-1}{m}
$$
and  $\alpha_{1} = j \bmod m$, 
\begin{align}\label{lj}
\ell(j)=-1-\frac{l\left(\alpha_{1}-j\right)}{m}-E_{j},
\end{align}
where 
$ E_{j}=\left\lfloor\frac{(j+n)l-1}{m}\right\rfloor$.
Therefore,
\begin{align*}
	(\ref{Kj1})=& \frac{1}{\pi^{n+1}} \sum_{\alpha_{1} \in \mathcal{R}_{j}} \sum_{\beta=\ell(j)}^{\infty}\left[\frac{\Gamma(\alpha_1+n+1)(\alpha_1+n)l}{\Gamma(\alpha_1+1)m}+\frac{\Gamma(\alpha_1+n+1)(\beta+1)}{\Gamma(\alpha_1+1)}\right] a^{\alpha_{1}} b^{\beta} \\
	=& \frac{l}{\pi^{n+1}m} \sum_{\alpha_{1} \in \mathcal{R}_{j}} \sum_{\beta=\ell(j)}^{\infty}\frac{\Gamma(\alpha_1+n+1)(\alpha_1+n)}{\Gamma(\alpha_1+1)} a^{\alpha_{1}} b^{\beta} \\
	&+ \frac{1}{\pi^{n+1}}\sum_{\alpha_{1} \in \mathcal{R}_{j}} \sum_{\beta=\ell(j)}^{\infty}\frac{\Gamma(\alpha_1+n+1)(\beta+1)}{\Gamma(\alpha_1+1)} a^{\alpha_{1}} b^{\beta} \\
	:=& \frac{l}{m \pi^{n+1}} I(j)+\frac{1}{\pi^{n+1}} J(j) .
\end{align*}
Next, we will calculate the sums $I(j)$ and $J(j)$ separately. Let $u:=a b^{-l / m}$, and note that both $a^m<b^l<1$ and $|u|<1$. 
\begin{align}\label{Ij}
	I(j)=&\sum_{\alpha_{1} \in \mathcal{R}_{j}}\frac{\Gamma(\alpha_1+n+1)(\alpha_1+n)}{\Gamma(\alpha_1+1)} a^{\alpha_{1}} \sum_{\beta=\ell(j)}^{\infty} b^{\beta}\notag\\ 
	=&\frac{1}{1-b} \cdot \sum_{\alpha_{1} \in \mathcal{R}_{j}}\frac{\Gamma(\alpha_1+n+1)(\alpha_1+n)}{\Gamma(\alpha_1+1)} a^{\alpha_{1}} b^{\ell(j)} \notag\\
	=&\frac{b^{\frac{l j}{m}-1-E_{j}}}{1-b} \cdot \sum_{\alpha_{1} \in \mathcal{R}_{j}}\frac{\Gamma(\alpha_1+n+1)(\alpha_1+n)}{\Gamma(\alpha_1+1)}  u^{\alpha_{1}} \notag\\
	=&\frac{b^{\frac{l j}{m}-1-E_{j}}}{1-b} \cdot \frac{\mathrm{d}^n}{\mathrm{d} u^n}\left(u \frac{\mathrm{d}}{\mathrm{d} u}\left(\frac{u^{j+n}}{1-u^{m}}\right)\right)\notag\\
	=&\frac{b^{\frac{l j}{m}-1-E_{j}}}{1-b} \cdot\left[n \cdot \frac{\mathrm{d}^n}{\mathrm{d} u^n}\left(\frac{u^{j+n}}{1-u^{m}}\right)+ u \cdot \frac{\mathrm{d}^{n+1}}{\mathrm{d} u^{n+1}}\left(\frac{u^{j+n}}{1-u^{m}}\right)\right].
\end{align}
The sum $J(j)$ is split into two pieces:
$$
\begin{aligned}
&~~J(j) \\
=&\sum_{\alpha_{1} \in \mathcal{R}_{j}} \frac{\Gamma(\alpha_1+n+1)}{\Gamma(\alpha_1+1)} a^{\alpha_{1}}  \sum_{\beta=\ell(j)}^{\infty}\left(\beta+1\right) b^{\beta} \\
=	&\sum_{\alpha_{1} \in \mathcal{R}_{j}} \frac{\Gamma(\alpha_1+n+1)}{\Gamma(\alpha_1+1)} a^{\alpha_{1}}  \sum_{\beta=\ell(j)}^{\infty} \frac{\mathrm{d}}{\mathrm{d} b} \left(b^{\beta+1} \right)\\
=	&\sum_{\alpha_{1} \in \mathcal{R}_{j}}\frac{\Gamma(\alpha_1+n+1)}{\Gamma(\alpha_1+1)} a^{\alpha_{1}}\left[\frac{b^{\ell(j)+1}}{(1-b)^{2}}+\frac{(\ell(j)+1) b^{\ell(j)}}{1-b}\right] \\
	=&\frac{b}{(1-b)^{2}} \sum_{\alpha_{1} \in \mathcal{R}_{j}}\frac{\Gamma(\alpha_1+n+1)}{\Gamma(\alpha_1+1)} a^{\alpha_{1}} b^{\ell(j)}+\frac{1}{1-b} \sum_{\alpha_{1} \in \mathcal{R}_{j}}\left(\ell(j)+1\right)\frac{\Gamma(\alpha_1+n+1)}{\Gamma(\alpha_1+1)} a^{\alpha_{1}} b^{\ell(j)} \\
:=	&J_{1}(j)+J_{2}(j) .
\end{aligned}
$$
For the first piece, it follows
\begin{align}\label{J1}
	J_{1}(j) =\frac{b^{\frac{l j}{m}-E_{j}}}{(1-b)^{2}} \sum_{\alpha_{1} \in \mathcal{R}_{j}}\frac{\Gamma(\alpha_1+n+1)}{\Gamma(\alpha_1+1)} u^{\alpha_{1}}
	&=\frac{b^{\frac{l j}{m}-E_{j}}}{(1-b)^{2}} \sum_{\alpha_{1} \in \mathcal{R}_{j}}
	\frac{\mathrm{d}^n}{\mathrm{d} u^n} \left(u^{\alpha_{1}+n}\right) \notag\\
	&=\frac{b^{\frac{l j}{m}-E_{j}}}{(1-b)^{2}} \cdot \frac{\mathrm{d}^n}{\mathrm{d} u^n}\left(\frac{u^{j+n}}{1-u^{m}}\right),
\end{align}
where $u=a b^{-l / m}$. For the second piece,
\begin{align}\label{J2}
J_{2}(j) 
    =&\frac{b^{\frac{l j}{m}-1-E_{j}}}{1-b} \sum_{\alpha_{1}\in \mathcal{R}_{j}}\left(\ell(j)+1\right)\frac{\Gamma(\alpha_1+n+1)}{\Gamma(\alpha_1+1)} u^{\alpha_{1}} \\
	=&\frac{b^{\frac{l j}{m}-1-E_{j}}}{1-b}\left[\left(\frac{l j}{m}-E_{j}\right) \sum_{\alpha_{1} \in \mathcal{R}_{j}}\frac{\Gamma(\alpha_1+n+1)}{\Gamma(\alpha_1+1)} u^{\alpha_{1}}-\frac{l}{m} \sum_{\alpha_{1} \in \mathcal{R}_{j}}\frac{\Gamma(\alpha_1+n+1)}{\Gamma(\alpha_1+1)} \alpha_{1} u^{\alpha_{1}}\right]  \notag\\
	=&\frac{b^{\frac{l j}{m}-1-E_{j}}}{1-b}\left[\left(\frac{l j}{m}-E_{j}\right) \cdot \frac{\mathrm{d}^n}{\mathrm{d} u^n}\left(\frac{u^{j+n}}{1-u^{m}}\right)-\frac{l}{m} \cdot u \cdot \frac{\mathrm{d}^{n+1}}{\mathrm{d} u^{n+1}}\left(\frac{u^{j+n}}{1-u^{m}}\right)\right]. \notag
\end{align}
 (\ref{Ij}) and (\ref{J2}) can be combined as
$$
I(j)+\frac{m}{l} J_{2}(j)=\frac{b^{\frac{l j}{m}-1-E_{j}}}{1-b}\left(j+n-\frac{m}{l} E_{j}\right) \cdot \frac{\mathrm{d}^n}{\mathrm{d} u^n}\left(\frac{u^{j+n}}{1-u^{m}}\right).
$$
Combining this with (\ref{J1}), we have
\begin{align}\label{Kj2}
	K_{j}((z,w),(z,w))
	=&\frac{l}{m \pi^{n+1}}\left[I(j)+\frac{m}{l} J_{2}(j)+\frac{m}{l} J_{1}(j)\right] \notag\\
	=&\frac{l}{m \pi^{n+1}} \cdot g_{j}(b) \cdot \frac{b^{\frac{l j}{m}-1-E_{j}}}{(1-b)^{2}} \cdot \frac{\mathrm{d}^n}{\mathrm{d} u^n}\left(\frac{u^{j+n}}{1-u^{m}}\right),
\end{align}
where $g_{j}(b):=j+n-\frac{m}{l} E_{j}+\left(\frac{m}{l}+\frac{m}{l} E_{j}-j-n\right) b$. Using Leibniz's rule,
\begin{align*}
    \frac{\mathrm{d}^n}{\mathrm{d}u^n}\left(\frac{u^{j+n}}{1-u^{m}}\right)
    =\sum_{k=0}^{n}\binom{n}{k}\left(\frac{1}{1-u^{m}}\right)^{(k)}(u^{j+n})^{(n-k)}
    =\frac{u^j\cdot Q(u^m)}{(1-u^m)^{n+1}},	
\end{align*}
where $Q$ is a polynomial of degree at most $n$.
Note that $u=a b^{-l / m}$, 
\begin{align}\label{Kj3}
	(\ref{Kj2}) 
	&=\frac{l}{m \pi^{n+1}} \cdot g_{j}(b) \cdot \frac{b^{\frac{l j}{m}-1-E_{j}}}{(1-b)^{2}} \cdot \frac{u^j}{(1-u^m)^{n+1}}\cdot Q(u^m) \notag\\
	&=\frac{l}{m \pi^{n+1}} \cdot g_{j}(b)\cdot Q(a^m b^{-l}) \frac{a^jb^{(n+1)l-1-E_{j}}}{(1-b)^{2}\left(b^{l}-a^{m}\right)^{n+1}}.
\end{align}
Polarization now gives the formula for $K_{j}((z,w),(s,t))$, substituting $a=z\cdot\overline{s}$ and
$b=w\cdot\overline{t}$  into equation (\ref{Kj3}). Finally, note that $\Omega^{n+1}_{m/l}$ is a bounded domain where $|z|^m < |w|^l<1$ and the estimates
$$
|g_{j}(b)|\lesssim 1, \quad
|Q(a^m b^{-l})|\lesssim 1,
$$
then the sub-Bergman kernel $K_{j}$ satisfies the estimate
$$
|K_{j}((z,w),(s,t))|\lesssim \frac{|b|^{\frac{lj}{m}+(n+1)l-E_j-1}}{|1-b|^{2}|b^l-a^{m}|^{n+1}}.
$$
\end{proof}

Recall that $E_{j}=\left\lfloor\frac{(j+n)l-1}{m}\right\rfloor,$ so
$$
\frac{(j+n)l-1}{m}-1<E_{j} \leq \frac{(j+n)l-1}{m} \quad \forall j \in\{0, \ldots, m-1\}.
$$
Then, Theorem
\ref{thm3} and (\ref{decompose}) yield the following estimate on the full Bergman kernel:
\begin{cor}\label{cor kernel estimate}
  The Bergman kernel of the domain $\Omega^{n+1}_{m/l}$ satisfies the estimate
\begin{align*}
	|B_{\frac{m}{l},n}((z,w),(s,t))|\lesssim \frac{|b|^{(n+1)l-1-\frac{nl-1}{m}}}{|1-b|^{2}|b^l-a^{m}|^{n+1}},
\end{align*}
where 
$a=z\cdot\overline{s}$  and
$ b=w\cdot\overline{t}$.
\end{cor}

\subsection{ Boundedness of operators on $\Omega^{n+1}_{m/l}$}
If $\Omega \subset \mathbb{C}^{n+1}$ is a domain and $K$ is an a.e. positive, measurable function on $\Omega \times \Omega$, let $\mathcal{K}$ denote the integral operator with kernel $K$:
$$
\mathcal{K}(f)(z,w)=\int_{\Omega} K((z, w),(s,t)) f(s,t) \mathrm{d} V(s,t).
$$
The fundamental result regarding Boundedness of operators is as follows:
\begin{thm}\label{pro1}
  If the kernel of the domain $\Omega^{n+1}_{m/l}$ satisfies the estimate
\begin{align}\label{KA}
  |K((z,w),(s,t))|\lesssim \frac{|b|^A}{|1-b|^{2}|b^l-a^{m}|^{n+1}},
\end{align}
  where 
  $a=z\cdot\overline{s},$ 
  $ b=w\cdot\overline{t},$ 
   then, $\mathcal{K}: L^{p}\left(\Omega^{n+1}_{m/l}\right) \longrightarrow$ $L^{p}\left(\Omega^{n+1}_{m/l}\right)$ boundedly if
\begin{align}\label{pA}
\frac{2 nl+2 m}{A m+2nl+2 m-(n+1)l m}<p<\frac{2 nl+2 m}{(n+1)l m-A m},
\end{align}
when both denominators in (\ref{pA}) are positive and $A m+2nl+2 m-(n+1)l m>$ $(n+1)l m-A m$.
\end{thm}

Some lemmas are needed before proving Theorem \ref{pro1}.

\begin{lem}[Schur's Lemma {\rm\cite{LD2}}]\label{Schur}
	Let $\Omega \subset \mathbb{C}^{n}$ be a domain, $K$ be an a.e. positive, measurable function on $\Omega \times \Omega$, and $\mathcal{K}$ be the integral operator with kernel $K .$ Suppose there exists a positive auxiliary function $h$ on $\Omega$, and numbers $0<a<b$ such that for all $\epsilon \in[a, b)$, the following estimates hold:
	$$
	\begin{aligned}
		\mathcal{K}\left(h^{-\epsilon}\right)(z) &:=\int_{\Omega} K(z, w) h(w)^{-\epsilon} d V(w) \lesssim h(z)^{-\epsilon} \\
		\mathcal{K}\left(h^{-\epsilon}\right)(w) &:=\int_{\Omega} K(z, w) h(z)^{-\epsilon} d V(z) \lesssim h(w)^{-\epsilon} .
	\end{aligned}
	$$
	Then, $\mathcal{K}$ is a bounded operator on $L^{p}(\Omega)$ for all $p \in\left(\frac{a+b}{b}, \frac{a+b}{a}\right)$.
\end{lem}

\begin{lem}{\rm\cite{LD2}}
	\label{estimate1}
	Let $D \subset \mathbb{C}$ be the unit disk, $\epsilon \in(0,1)$ and $\beta \in(-\infty, 2)$. Then, for $z \in D$,
	$$
	\mathcal{I}_{\epsilon, \beta}(z):=\int_{D} \frac{\left(1-|w|^{2}\right)^{-\epsilon}}{|1-z \bar{w}|^{2}}|w|^{-\beta} \mathrm{d} V(w) \lesssim\left(1-|z|^{2}\right)^{-\epsilon}
	$$
	with a constant independent of $z$.
\end{lem}

\begin{lem}\label{estimate2}{\rm\cite{FU}}
	Let $D_n \subset \mathbb{C}^n$ be the unit ball, $k\in \mathbb{Z}^+$, $\epsilon \in(0,1)$ and $\Delta \in \mathbb{C}^n, |\Delta|<1 $. Then
	\begin{align}\label{lem9}
		\int_{D_n} \frac{\left(1-|\eta|^{2k}\right)^{-\epsilon}}{|1-(\eta \cdot \overline{\Delta})^k|^{n+1}} \mathrm{d} V(\eta) \approx \left(1-|\Delta|^{2k}\right)^{-\epsilon}.
	\end{align}
\end{lem}

$~$

\begin{proof}[\textbf{Proof of Theorem  $\mathbf{\ref{pro1}}$}]
	Let $h(z,w):= (|w|^{2l}-|z|^{2m})(1-|w|^2)$. This function (essentially) measures the distance of $(z,w) \in \Omega^{n+1}_{m/l}$ to $b(\Omega^{n+1}_{m/l})$.
	We will prove that for all $\epsilon \in [\frac{(n+1)}{2}-\frac{A}{2l},\frac{n}{m}+\frac{1}{l}+\frac{A}{2l}-\frac{n+1}{2})$, and any $(z,w) \in \Omega^{n+1}_{m/l}$,
	\begin{align}\label{th4.2}
		|\mathcal{K}|\left(h^{-\epsilon}\right)(z,w) :=\int_{\Omega^{n+1}_{m/l}} |K((z, w),(s,t))| h(s,t)^{-\epsilon} d V(s,t) \lesssim h(z,w)^{-\epsilon}.
	\end{align}
	From estimate (\ref{KA}), we see that
	\begin{align*}
		|\mathcal{K}|\left(h^{-\epsilon}\right)(z,w)
		&\lesssim  \int_{0<|t|<1}\int_{|s|^{\frac{m}{l}}<|t|}\frac{|w\cdot\overline{t}|^{A}(|t|^{2l}-|s|^{2m})^{-\epsilon}(1-|t|^{2})^{-\epsilon}}{|1-w\cdot\overline{t}|^{2}|w^l\cdot\overline{t^l}-(z\cdot\overline{s})^{m}|^{n+1}} dV(s)dV(t).
	\end{align*}
	Let
	\begin{align*}
		t=\rho e^{i\varphi},\quad
		w=|w| e^{i\varphi_0}, \quad
		s=(r_1 e^{i\theta_1}, \cdots, r_n e^{i\theta_n}), \quad
		z=(|z_1| e^{i\tilde{\theta}_1}, \cdots, |z_n| e^{i\tilde{\theta}_n}),
	\end{align*}
	then
		\begin{align*}
			&~~|\mathcal{K}|\left(h^{-\epsilon}\right)(z,w)\\
			\lesssim & \int_{0}^1\int_{|r|<\rho^{\frac lm}}(\prod_{k=1}^{n+1}\int_{-\pi}^{\pi})
			\frac{|w|^{A}\rho^{A+1}(\rho^{2l}-|r|^{2m})^{-\epsilon}(1-\rho^{2})^{-\epsilon}r_1\cdots r_nd\varphi d\theta_1\cdots d\theta_ndrd\rho}
			{|1-|w|\rho e^{i(\varphi_0-\varphi)}|^{2}||w|^l\rho^l e^{il(\varphi_0-\varphi)}-
			\Upsilon_1
			|^{n+1}}\\
			=&\int_{0}^1\int_{|r|<\rho^{\frac lm}}(\prod_{k=1}^{n+1}\int_{-\pi}^{\pi})
			\frac{|w|^{A}\rho^{A+1}(\rho^{2l}-|r|^{2m})^{-\epsilon}(1-\rho^2)^{-\epsilon}r_1\cdots r_nd\varphi d\theta_1\cdots d\theta_ndrd\rho}
			{|1-|w|\rho e^{i\varphi}|^{2}||w|^l\rho^l-(|z_1|r_1 e^{i\theta_1}+\cdots+|z_n|r_n e^{i\theta_n})^{m}|^{n+1}},
		\end{align*}
where 
$\Upsilon_1=\left( |z_1|r_1 e^{i(\tilde{\theta}_1-\theta_1)}+\cdots+|z_n|r_n e^{i(\tilde{\theta}_n-\theta_n)}\right) ^{m}$.
In the last equation, we used the periodicity of the $\theta_1,\cdots, \theta_n$ and $\varphi$ integrals. Next, we first consider the following integral
	\begin{align}\label{th4.3}
		&~~\int_{|r|<\rho^{\frac lm}}
		(\prod_{k=1}^{n}\int_{-\pi}^{\pi})
		\frac{(\rho^{2l}-|r|^{2m})^{-\epsilon}r_1\cdots r_n d\theta_1\cdots d\theta_ndr}
		{||w|^l\rho^l-(|z_1|r_1 e^{i\theta_1}+\cdots+|z_n|r_n e^{i\theta_n})^{m}|^{n+1}}\\
		=&\int_{|r|<\rho^{\frac lm}}(\prod_{k=1}^{n}\int_{-\pi}^{\pi})
		\frac{\rho^{-2l\epsilon}(1-|r\rho^{-l/m}|^{2m})^{-\epsilon}r_1\cdots r_n d\theta_1\cdots d\theta_ndr}
		{(|w|\rho)^{(n+1)l}
		|1 -\Upsilon_2
		|^{n+1}},\notag
	\end{align}
	where $\Upsilon_2=\left( |z_1|r_1 (|w|\rho)^{-l/m} e^{i\theta_1}+\cdots+|z_n|r_n (|w|\rho)^{-l/m} e^{i\theta_n}\right) ^{m}$.
	{Make} the substitution $\tilde{r}_k=r_k\rho^{-l/m}$ $ (k=1,\cdots,n),$ and let
	$$\eta=(\tilde{r}_1e^{i\theta_1},\cdots, \tilde{r}_ne^{i\theta_n}),\quad
	\Delta=(|z_1||w|^{-l/m},\cdots, |z_n||w|^{-l/m}),$$ then
	\begin{align*}
		(\ref{th4.3})
		&=\int_{|\tilde{r}|<1}(\prod_{k=1}^{n}\int_{-\pi}^{\pi})
		\frac{\rho^{\frac{2nl}m-2l\epsilon-(n+1)l}(1-|\tilde{r}|^{2m})^{-\epsilon}\tilde{r}_1\cdots\tilde{r_n} d\theta_1\cdots d\theta_nd\tilde{r}}
		{|w|^{(n+1)l}|1 -(|z_1| |w|^{-l/m} \tilde{r}_1 e^{i\theta_1}+\cdots+|z_n| |w|^{-l/m} \tilde{r}_n e^{i\theta_n})^{m}|^{n+1}}\\
		&=\frac{\rho^{\frac{2nl}m-2l\epsilon-(n+1)l}}{|w|^{(n+1)l}}\int_{D_n} \frac{\left(1-|\eta|^{2m}\right)^{-\epsilon}}{|1-(\eta \cdot \overline{\Delta})^m|^{n+1}} \mathrm{d} V(\eta).
	\end{align*}
	From Lemma \ref{estimate2}, we have
	\begin{align*}
		(\ref{th4.3})\approx \frac{\rho^{\frac{2nl}m-2l\epsilon-(n+1)l}}{|w|^{(n+1)l}}\left(1-|\Delta|^{2m}\right)^{-\epsilon}
		=\rho^{\frac{2nl}m-2l\epsilon-(n+1)l}|w|^{2l\epsilon-(n+1)l}\left(|w|^{2l}-|z|^{2m}\right)^{-\epsilon}.
	\end{align*}
	This means that
	\begin{align*}
		&~~|\mathcal{K}|\left(h^{-\epsilon}\right)(z,w)\\
		\lesssim& |w|^{2l\epsilon+A -(n+1)l}\left(|w|^{2l}-|z|^{2m}\right)^{-\epsilon}
		\int_{0}^1\int_{-\pi}^{\pi}
		\frac{\rho^{\frac{2nl}m-2l\epsilon-(n+1)l+A+1}(1-\rho^2)^{-\epsilon}d\varphi d\rho}
		{|1-|w|\rho e^{i\varphi}|^{2}}\\
		=&|w|^{2l\epsilon+A -(n+1)l}\left(|w|^{2l}-|z|^{2m}\right)^{-\epsilon}
		\int_{D}
		\frac{|z|^{\frac{2nl}m-2l\epsilon-(n+1)l+A}(1-|z|^2)^{-\epsilon}dV(z)}
		{\left|1-|w|z\right|^{2}}.
	\end{align*}
	From Lemma \ref{estimate1},
	\begin{align*}
		|\mathcal{K}|\left(h^{-\epsilon}\right)(z,w)
		\lesssim |w|^{2l\epsilon+A -(n+1)l}\left(|w|^{2l}-|z|^{2m}\right)^{-\epsilon}\left(1-|w|^{2}\right)^{-\epsilon},
	\end{align*}
	when $\frac{2nl}m-2l\epsilon-(n+1)l+A > -2, i.e.,$ $ \epsilon<\frac{n}{m}+\frac{1}{l}+\frac{A}{2l}-\frac{n+1}{2}$. Then,
	\begin{align*}
		|\mathcal{K}|\left(h^{-\epsilon}\right)(z,w)
		\lesssim\left(|w|^{2l}-|z|^{2m}\right)^{-\epsilon}\left(1-|w|^{2}\right)^{-\epsilon}=h(z,w)^{-\epsilon},
	\end{align*}
	when $2l\epsilon+A -(n+1)l \geq0, i.e.$, $\frac{(n+1)}{2}-\frac{A}{2l}\leq \epsilon$.  This completes the proof of (\ref{th4.2}). Finally, combining (\ref{th4.2}) and Schur's Lemma (Lemma \ref{Schur}) yields that the operator $|\mathcal{K}|$ is bounded from $L^p(\Omega^{n+1}_{m/l})$ to $L^p(\Omega^{n+1}_{m/l})$ when
    $$
		\frac{2 nl+2 m}{A m+2nl+2 m-(n+1)l m}<p<\frac{2 nl+2 m}{(n+1)l m-A m}.
    $$
    Note that because of the conjugate symmetry of $K$, it is sufficient to establish just one of the estimates to apply Lemma \ref{Schur}.
    A fortiori, $\mathcal{K}$ is bounded from $L^p(\Omega^{n+1}_{m/l})$ to $L^p(\Omega^{n+1}_{m/l})$ for $p$ in the same range.  This completes the proof.
    	
\end{proof}

\begin{thm}
	\label{necessity}
   	The Bergman projection $\mathbf{P}_{m / l}$ is a bounded operator on $L^p(\Omega^{n+1}_{m/l})$ for   all $p \in\left(\frac{2 m+2 nl}{m+nl+1}, \frac{2 m+2 nl}{m+nl-1}\right)$.
\end{thm}

\begin{proof}
	It comes immediately from Corollary \ref{cor kernel estimate} and Theorem $\ref{pro1}$ by taking $A=(n+1)l-1-\frac{nl-1}{m}$.
\end{proof}

\begin{rem}
The proof via Schur's Lemma actually yields the stronger result that the operator $\mathbf{P}^{+}_{m/l}$, defined by replacing the Bergman kernel with its absolute value, is also bounded on the same $L^p(\Omega^{n+1}_{m/l})$ spaces.
\end{rem}

\section{The Rational Case: $L^p$ Unboundedness}\label{Qnonboundedness}
We shall show that $\mathbf{P}_{m/l}$ fails to be bounded on $L^p(\Omega^{n+1}_{m/l})$ for $p \notin\left(\frac{2 m+2 nl}{m+nl+1}, \frac{2 m+2 nl}{m+nl-1}\right)$ by exhibiting a single function $f \in L^{\infty}(\Omega^{n+1}_{m/l})$ such that $\mathbf{P}_{m / l} f \notin L^p(\Omega^{n+1}_{m/l})$.

\begin{lem}\label{function}
Suppose $(\eta, \nu) \in \mathbb{N}^n \times \mathbb{Z}$ is an index such that both $z^\eta w^\nu$ and $z^\eta w^{-\nu}$ are elements of $A^2\left(\Omega_\gamma^{n+1}\right)$. If we let $f(z, w):=z^\eta \bar{w}^\nu$, there exists a positive constant $C$ such that

$$
\mathbf{P}_\gamma(f)(z, w)=C z^\eta w^{-\nu} .
$$
\end{lem}

\begin{proof}
	From Lemma \ref{A2}, the Bergman kernel admits the expansion
	\begin{align*}
			B_{\gamma,n}((z,w),(s,t))
			=\sum\limits_{|\alpha| + \gamma (\beta + 1) > -n}\frac{z^{\alpha}w^{\beta}(\overline{s^\alpha t^\beta})}{\|z^\alpha w^\beta\|^2_{L^2(\Omega^{n+1}_\gamma)}}.
	\end{align*}
	It follows that
	\begin{align*}
		\mathbf{P}_{\gamma}(f)(z,w) =&\int_{\Omega^{n+1}_{\gamma}} B_{\gamma,n}((z, w),(s,t)) f(s,t) d V(s,t)\\
		=&\int_{\Omega^{n+1}_{\gamma}} \sum\limits_{|\alpha| + \gamma (\beta + 1) > -n}\frac{z^{\alpha}w^{\beta}(\overline{s^\alpha t^\beta})}{\|z^\alpha w^\beta\|^2_{L^2(\Omega^{n+1}_\gamma)}}s^{\eta}\bar{t}^{\nu} d V(s,t) \\
		=&\sum\limits_{|\alpha| + \gamma (\beta + 1) > -n} 
		\frac{z^\alpha w^\beta}{\|z^\alpha w^\beta\|^2_{L^2(\Omega^{n+1}_\gamma)}} \int_{{\Omega^{n+1}_{\gamma}}} \overline{s^\alpha t^\beta}  s^\eta \overline{t}^\nu d V(s,t).
	\end{align*}
	Since $\Omega^{n+1}_\gamma$ is a Reinhardt domain, we apply polar coordinate transformations to each component separately. Let
	\begin{align*}
		t=\rho e^{i\varphi},\quad
		s=(r_1 e^{i\theta_1}, \cdots, r_n e^{i\theta_n}), \quad
		r=(r_1, \cdots, r_n),\quad
		\theta=(\theta_1, \cdots, \theta_n),
	\end{align*}
	then 
		\begin{align*}
		&\int_{{\Omega^{n+1}_{\gamma}}} \overline{s^\alpha t^\beta}  s^\eta \overline{t}^\nu d V(s,t)\\
		=&\int_{0}^1\int_{|r|^{\gamma}<\rho}(\prod_{k=1}^{n+1}\int_{-\pi}^{\pi})r^{\alpha+\eta+\mathbbm{1}_n} \rho^{\beta+\nu+1} e^{i\theta_1(\eta_1-\alpha_1)} \cdots e^{i\theta_n(\eta_n-\alpha_n)}   e^{-i\varphi(\beta+\nu)} d\varphi d\theta_1\cdots d\theta_ndrd\rho,
	\end{align*}
where $\mathbbm{1}_n=(1,1,\cdots,1)$. The integrations vanish unless $(\alpha,\beta) = (\eta,-\nu)$. Thus,
 $
\mathbf{P}_\gamma(f)(z, w)=C z^\eta w^{-\nu},
$ where $C$ is a constant.
\end{proof}

\begin{cor}\label{pro2}
	If both $(\eta_1, \eta_2)$ and $(\eta_1, -\eta_2)$ belong to $\Lambda_{\gamma,n}$, and $f(z,w):=z_1^{\eta_1}\bar{w}^{\eta_2}.$ Then there exists a
	constant $C$ such that
	\begin{align*}
		\mathbf{P}_{\gamma}(f)(z,w)=Cz_1^{\eta_1}w^{-\eta_2}.
	\end{align*}	
\end{cor}
\begin{proof}
It comes immediately from  Lemma $\ref{function}$ by taking $(\eta, \nu)=\left(\eta_1, 0, \cdots, 0, \eta_2\right) \in \mathbb{N}^n \times \mathbb{Z}$.
\end{proof}

Recall that $\mathcal{K}_j$ denotes the sub-Bergman projection defined in Section \ref{Qboundedness}. For $\gamma \in \mathbb{Q}^+$, an analogous result for the subspaces $\mathcal{S}_{j}$ can be established using the same line of reasoning.

\begin{lem} \label{pro3}
	If both $\left(\eta_{1}, \eta_{2}\right)$ and $\left(\eta_{1},-\eta_{2}\right)$ belong to $\mathcal{G}_{j}$ for some $j \in$ $\{0,1, \ldots, m-1\}$, and $f(z,w):=z_1^{\eta_1}\bar{w}^{\eta_2},$ then there exists a constant $C$ such that
$$
\mathcal{K}_{l}(f)(z,w)= \begin{cases}Cz_1^{\eta_1}w^{-\eta_2}, & l=j \\ 0, & l \neq j\end{cases}
$$
for all $l \in\{0,1, \ldots, m-1\}$,where $\mathcal{G}_{j}=\left\{(\alpha_1,\beta) \in \Lambda_{\frac{m}{l},n}: \alpha_{1} = j \bmod m\right\}$.
\end{lem} 

\begin{thm}\label{pro4} 
	For each $j \in\{0,1, \ldots, m-1\}$, the sub-Bergman projection $\mathcal{K}_{j}$ does not map $L^{\infty}(\Omega^{n+1}_{m/l})$ to $L^{p}(\Omega^{n+1}_{m/l})$ for any $p \geq \frac{2 m+2 nl}{m+m E_{j}-l j}$, where  
	$ E_{j}=\left\lfloor\frac{(j+n)l-1}{m}\right\rfloor$.
\end{thm}
\begin{proof}
 Fix $j$, and take $\eta_{1}=j+k m$ for some $k \in \mathbb{Z}^{+} \cup\{0\}$. Let $\eta_{2}=\ell\left(\eta_{1}\right)$, and note that (\ref{lj}) says that
$$
\eta_{2}=-1-l k-E_{j}<0.
$$
Thus, $\left(\eta_{1}, \eta_{2}\right),\left(\eta_{1},-\eta_{2}\right) \in \mathcal{G}_{j}$. Let $f(z,w):=z_1^{\eta_1}/\bar{w}^{\eta_2}$; clearly $f \in L^{\infty}(\Omega^{n+1}_{m/l})$. Note Corollary \ref{pro2} says that $\mathcal{K}_{j} (f)(z,w)=Cz_1^{\eta_1}w^{\eta_2}$. Thus, 
\begin{align*}
	&~~\int_{\Omega^{n+1}_{m/l}} | z_1|^{{\eta_1}p}\left|w\right|^{{\eta_2}p} dV(z,w)\\
	=&\int_{0<|w|<1} |w|^{{\eta_2}p}\int_{0}^{|w|^{\frac{l}{m}}}\int_{|\xi|=1}r^{{\eta_1}p+2n-1}| \xi_1|^{{\eta_1}p} d\sigma(\xi)drdV(w) \\
	=&\frac{1}{2n+{\eta_1}p}\int_{|\xi|=1}|\xi_1|^{{\eta_1}p} d\sigma(\xi)\int_{0<|w|<1}|w|^{{\eta_2}p+\frac{l}{m}(2n+{\eta_1}p)}dV(w)\\
	=&\frac{2\pi }{2n+{\eta_1}p}\int_{|\xi|=1}|\xi_1|^{{\eta_1}p} d\sigma(\xi)\int_{0}^1r^{{\eta_2}p+\frac{l}{m}(2n+{\eta_1}p)+1}dr.
\end{align*}
This integral diverges if and only if
\begin{align}\label{prange}
{\eta_2}p+\frac{l}{m}(2n+{\eta_1}p)+1 \leq-1.
\end{align}
Substituting $\eta_{1}=j+k m$ and $\eta_{2}=-1-l k-E_{j}$, (\ref{prange}) becomes
\begin{align}
-p\left(m+m E_{j}-l j\right) \leq-2 nl-2 m .
\end{align}
Due to $E_{j}=\left\lfloor\frac{(j+n)l-1}{m}\right\rfloor$,
$$
\begin{aligned}
	m+m E_{j}-l j &>m+m\left( \frac{(j+n)l-1}{m}-1\right) -l j \\
	&=nl-1 \geq 0.
\end{aligned}
$$
Therefore,  $p \geq \frac{2 m+2 nl}{m+m E_{j}-l j}$.
\end{proof}

\begin{thm}\label{pro4.4}  
	For $p \geq \frac{2 m+2 nl}{m+nl-1}$,
	$\mathbf{P}_{m / l}$ fails to map $L^{\infty}(\Omega^{n+1}_{m/l})$ to $L^{p}(\Omega^{n+1}_{m/l})$.
\end{thm} 

\begin{proof}
	As $m$ and $l$ are relatively prime, according to elementary number theory, there is a unique $x \in \{n, \ldots, n+m-1\}$ such that 
	$$ 
	l x=1 \quad \text { mod}\  m.  
	$$
	Setting $j_0 = x- n$, thus
	\begin{align}\label{Ej0}
		E_{j_{0}}=\frac{(j_0+n)l-1}{m}
	\end{align}
	and
	\begin{align}\label{ell}
		\ell\left(j_{0}\right) =-1-\frac{l\left(j_{0}-j_{0}\right)}{m}-E_{j_{0}}=-\frac{l}{m} j_{0}-1+\frac{1-nl}{m}.
	\end{align}
    Thus, $\left(j_0, \ell(j_{0})\right),\left(j_0, -\ell(j_{0})\right) \in \mathcal{G}_{j_0}$.
	Theorem \ref{pro4} says that $\mathcal{K}_{j_{0}}$ does not map the bounded function $g(z,w)=z_{1}^{j_{0}} /\bar{w}^{\ell\left(j_{0}\right)}$ to $L^{p}(\Omega^{n+1}_{m/l})$ for $p \geq \frac{2 m+2 nl}{m+nl-1}$. 
	On the other hand, Lemma \ref{pro3} says that $\mathcal{K}_{j}(g)=0$ for all $j \neq j_0$. Consequently, equation (\ref{decompose}) directly yields the desired conclusion.
\end{proof}
For $p < 2$, consider a fundamental implication of the Bergman projection's self-adjoint property:
\begin{lem}{\rm\cite{LD2}}\label{LD}
	Let $\Omega$ be a bounded domain and  $p > 1$. If $\mathbf{P}$ maps $L^p(\Omega)$ to $L^p(\Omega)$
	boundedly, then it also maps $L^q(\Omega)$ to $L^q(\Omega)$ boundedly, where $\frac1p+\frac1q=1$.
\end{lem}
 \begin{thm}\label{Sufficiency}  
	$\mathbf{P}_{m / l}$ is not a bounded operator on $L^{p}(\Omega^{n+1}_{m/l})$ for $p \notin\left(\frac{2 m+2 nl}{m+nl+1}, \frac{2 m+2 nl}{m+nl-1}\right)$.
\end{thm} 
\begin{proof}
	Theorem \ref{pro4.4} and Lemma \ref{LD} give   Theorem  \ref{Sufficiency}.
\end{proof}


\section{The Irrational Case: Degenerate $L^p$ Mapping}\label{Irrational}

\begin{lem}[Theorem 193, Page 164 \rm\cite{GH}]\label{Dirichlet}
	If $\gamma$ is irrational, $n \in \mathbb{Z}^+,$ there exists a sequence of rational numbers $\left\{\frac{m_{k}}{l_{k}}\right\}$, with $\frac{m_{k}}{l_{k}} \rightarrow \gamma$, such that
	$$
	\left|\frac{l_{k}}{m_{k}}-\frac{1}{\gamma}\right|<\frac{1}{m_{k}^{2}\sqrt{5}} .
	$$
\end{lem}

Next, we use the above lemma to prove (2) in Theorem \ref{thm Bergman projection}.
\begin{proof}[\textbf{Proof of Theorem  $\mathbf{\ref{thm Bergman projection}}$ (2)}]
For $p>2$, we will exhibit an $f \in L^{\infty}\left(\Omega^{n+1}_{\gamma}\right)$ such that $\mathbf{P}_{\gamma}(f) \notin$ $L^{p}\left(\Omega^{n+1}_{\gamma}\right)$.
Let $\left\{\frac{m_{k}}{l_{k}}\right\}$ be a sequence of rational numbers given by Lemma \ref{Dirichlet} and $m_k>n-1$. 
From (\ref{ell}), there exists a unique $\eta=\left(\eta_{1}, \eta_{2}\right) \in \Lambda_{m_{k} / l_{k},n}$ with $0 \leq$ $\eta_{1} \leq m_{k}-1$ such that
\begin{align}\label{eta2}
	\eta_{2}=\frac{1-l_{k} \eta_{1}-nl_{k}-m_{k}}{m_{k}} \in \mathbb{Z}.
\end{align}
Assume for the moment that this multi-index $\eta \in \Lambda_{\gamma,n}$. We will explain at the end of the proof that this is always the case.

Let $f_{k}(z,w):=z_{1}^{\eta_{1}}/ \bar{w}^{\eta_{2}}$; as $\eta_{2}<0, f_{k} \in L^{\infty}\left(\Omega^{n+1}_{\gamma}\right)$. Since we are assuming $\eta \in \Lambda_{\gamma,n}$, Lemma \ref{pro2} implies $\mathbf{P}_{\gamma}\left(f_{k}\right)(z,w) \approx z_{1}^{\eta_{1}} w^{\eta_{2}}$. It follows that
$$
\begin{aligned}
	\left\|\mathbf{P}_{\gamma}\left(f_{k}\right)\right\|_{L^{p}\left(\Omega^{n+1}_{\gamma}\right)}^{p} 
	& \approx \int_{\Omega^{n+1}_{\gamma}}|z_{1}|^{\eta_{1} p} |w|^{\eta_{2} p} \mathrm{d} V(z,w)\\	&=\int_{0<|w|<1} |w|^{{\eta_2}p}\int_{0}^{|w|^{\frac{1}{\gamma}}}\int_{|\xi|=1}r^{{\eta_1}p+2n-1}| \xi_1|^{{\eta_1}p} d\sigma(\xi)drdV(w) \\
	&\approx \int_{0<|w|<1}|w|^{{\eta_2}p+\frac{1}{\gamma}(2n+{\eta_1}p)}dV(w)\\
	&\approx \int_{0}^1r^{{\eta_2}p+\frac{1}{\gamma}(2n+{\eta_1}p)+1}dr.
\end{aligned}
$$
This diverges if the exponent is ${\eta_2}p+\frac{1}{\gamma}(2n+{\eta_1}p)+1\leq-1$. Since $\eta_{2}=\frac{1-l_{k} \eta_{1}-nl_{k}-m_{k}}{m_{k}} $, this means 
\begin{align}\label{Degenerate}
p\left(1+\frac{nl_{k}-1}{m_{k}}+\eta_{1}\left(\frac{l_{k}}{m_{k}}-\frac{1}{\gamma}\right)\right) \geq 2+\frac{2n}{\gamma}.
\end{align}
Consider the left hand side of (\ref{Degenerate}). Since $0 \leq \eta_{1} \leq m_{k}-1$,
$$
\eta_{1}\left|\frac{l_{k}}{m_{k}}-\frac{1}{\gamma}\right|
<\frac{m_{k}-1}{m_{k}^{2}\sqrt{5}}
<\frac{1}{m_{k}},
$$
by Lemma \ref{Dirichlet}. Thus
$$
\begin{aligned}
	p\left(1+\frac{nl_{k}-1}{m_{k}}+\eta_{1}\left(\frac{l_{k}}{m_{k}}-\frac{1}{\gamma}\right)\right) & \geq p\left(1+\frac{nl_{k}-1}{m_{k}}-\eta_{1}\left|\frac{l_{k}}{m_{k}}-\frac{1}{\gamma}\right|\right) \\
	&>p\left(1+\frac{nl_{k}-2}{m_{k}}\right).
\end{aligned}
$$
However since $p>2$, we can always choose $k$ large enough so that
$$
p\left(1+\frac{nl_{k}-2}{m_{k}}\right)>2+\frac{2n}{\gamma}.
$$
Thus, (\ref{Degenerate}) is satisfied for such $k$, which shows $\mathbf{P}_{\gamma}\left(f_{k}\right) \notin L^{p}\left(\Omega^{n+1}_{\gamma}\right)$.
We now show that the unique multi-index $\eta=\left(\eta_{1}, \eta_{2}\right) \in \Lambda_{m_{k} / l_{k},n}$ with $0 \leq \eta_{1} \leq$ $m_{k}-1$ and $\eta_{2}$ given by $(\ref{eta2})$ is necessarily in $\Lambda_{\gamma,n}$. We leave off the subscript $k$ in what follows.

Again, the rational approximation $\left|\frac{l}{m}-\frac{1}{\gamma}\right|<\frac{1}{\sqrt{5}m^{2}}$ is essential. If $\frac{m}{l}>\gamma$, then $A^{2}\left(\Omega^{n+1}_{m / l}\right) \subset A^{2}\left(\Omega^{n+1}_{\gamma}\right)$, so automatically, $\eta \in \Lambda_{\gamma,n}$. Suppose instead that $\frac{m}{l}<\gamma$. (\ref{Lambda1}) implies that $\eta=(\eta_1, \eta_2) \in \Lambda_{\gamma,n}$ if and only if $\eta_{1} \geq 0$ and $\eta_{2} > g\left(\eta_{1}\right)$, where
$$
g\left(\eta_{1}\right):=-\frac{\eta_{1}}{\gamma}-\frac{n}{\gamma}-1.
$$
However,  since $\frac{m}{l} \in \mathbb{Q}^{+}$, a multi-index $\eta=(\eta_1, \eta_2) \in \Lambda_{m / l, n}$ if and only if both $\eta_{1} \geq 0$ and $\eta_{2} > h\left(\eta_{1}\right)$, where
$$
h\left(\eta_{1}\right):=-\frac{l}{m} \eta_{1}+\frac{1-nl}{m}-1.
$$
Now under the requirement of $m>n-1$, for $0 \leq \eta_{1} \leq m-1$,
$$
\begin{aligned}
	h\left(\eta_{1}\right)-g\left(\eta_{1}\right) &=\frac{1}{m}-\left(\eta_{1}+n\right)\left(\frac{l}{m}-\frac{1}{\gamma}\right) \\
	& \geq \frac{1}{m}-(m+n-1)\frac{1}{\sqrt{5}m^2} \\
	&>0 .
\end{aligned}
$$
From this, it follows that $\eta=\left(\eta_{1}, \eta_{2}\right) \in \Lambda_{\gamma,n}$.
Since $p>2$ was arbitrary, the above shows that $\mathbf{P}_{\gamma}$ is not $L^{p}$ bounded for any $p>2$. Lemma \ref{LD} now shows that $\mathbf{P}_{\gamma}$ is not $L^{p}$ bounded for any $1<p<2$, which completes the proof.  
\end{proof}

\section*{Acknowledgments}
 The authors are deeply grateful to the anonymous referee for their
valuable comments and careful review. Their thoughtful feedback has been instrumental
in strengthening the paper.
\subsection*{Funding}
The project is supported by the National Natural Science Foundation of China (Grant no. 12571085), National Key Research and Development Project of China (Grant no. 2021YFA1002600) and Natural Science Foundation of Beijing Municipal (No.1252005).

\subsection*{Data Availability} Our manuscript has no associated data.

\section*{Declarations}
\subsection*{Conflicts of Interest} No potential conflict of interest was reported by author(s).

\vskip 0.5cm{
\end{document}